\documentclass[11pt]{article}
\usepackage{hyperref}
\usepackage{float}
\usepackage{amsf  onts,amsmath,amssymb}
\usepackage{xcolor}
\usepackage{graphicx}
\usepackage{setspace}
\def\bR{{\mathbf{R}}}

\newcommand{\ese}{\end{eqnarray*}}
\newcommand{\bse}{\begin{eqnarray*}}

\newtheorem{proposition}{Proposition}[section]

\newtheorem{definition}{Definition}[section]
\newtheorem{example}{Example}[section]

\newcommand{\cP}{{\cal P}}
\def\qed{\ \hfill$\square$\par\smallskip}

\def\Conv{\hbox{\rm Conv}}

\newcommand{\I}{{\cal I}}

\newcommand{\be}{\begin{eqnarray}}
\newcommand{\ee}[1]{\label{#1}\end{eqnarray}}

\newcommand{\rf}[1]{~(\ref{#1})}

\def\pen{{\hbox{\scriptsize\rm pen}}}

\def\cS{{\cal S}}

\def\cI{{\cal I}}

\def\bR{{\mathbf{R}}}

\def\bE{{\mathbf{E}}}

\newcommand{\aic}[2]{{\color{cyan}~#2}}

\title{Constant Depth Decision Rules for multistage optimization under uncertainty}

\author{Vincent Guigues
 \thanks{Corresponding author, School of Applied Mathematics FGV/EMAp, 22 250-900 Rio de Janeiro, Brazil
 {\tt vincent.guigues@fgv.br}}
 \and
Anatoli Juditsky
\thanks{LJK, Universit\'e Grenoble Alpes, 700 Avenue Centrale 38041 Domaine Universitaire
de Saint-Martin-d'H\`{e}res, France,
{\tt anatoli.juditsky@imag.fr}}
\and Arkadi Nemirovski
\thanks{Georgia Institute of Technology, Atlanta, Georgia
30332, USA, {\tt nemirovs@isye.gatech.edu}}}
\date{}

\begin{document}

\maketitle

\begin{abstract} In this paper, we introduce a new class of decision rules, referred to as {\em Constant Depth Decision Rules} (CDDRs),  for multistage optimization under linear constraints with uncertainty-affected
right-hand sides.
 We consider two uncertainty classes: discrete uncertainties which can take at each stage at most a fixed number $d$ of different values, and polytopic uncertainties which, at each stage, are elements of a convex hull of at most $d$ points. Given the {\em depth} $\mu$ of the decision rule, the decision at stage $t$
is expressed as the sum of $t$ functions of $\mu$ consecutive values of the underlying uncertain parameters. These functions are arbitrary in the case  of discrete
uncertainties and are poly-affine in the case of polytopic
uncertainties. For these uncertainty classes, we show that when the uncertain right-hand sides of the constraints of the multistage problem are of the same additive structure as the decision rules, these constraints can be reformulated as a system of linear inequality constraints where the numbers of variables and constraints is $O(1)(n+m)d^\mu N^2$ with $n$ the maximal dimension of control variables, $m$ the maximal number of inequality constraints at each stage, and $N$ the number of stages.
\par As an illustration, we discuss an application of the proposed approach to a Multistage Stochastic Program arising in the problem of hydro-thermal production planning with interstage dependent inflows. For  problems with a small number of stages, we present the results of a numerical study in which optimal CDDRs show similar performance, in terms of optimization objective, to that of Stochastic Dual Dynamic Programming (SDDP) policies, often at much smaller computational cost.
\end{abstract}



\section{Introduction}
Multistage optimization problems under uncertainty arise in
many real-life applications in finance and engineering, see for instance \cite{birge-louv-book, shadenrbook} and references therein, but are challenging to solve.
Such problems with stochastic uncertainties---Multistage Stochastic Programs (MSPs)---are typically replaced with their discretized (scenario) formulations \cite{shadenrbook}. Because the number of variables
in the deterministic equivalent of such approximations increases very quickly with the number $N$ of stages these problems are computationally intractable \cite{shapnemcompl} in general.
``Practical'' solution methods for convex MSPs often use decomposition techniques based on Dynamic Programming with two popular methods being Approximate Dynamic
Programming and Stochastic Dual Dynamic Programming (SDDP), introduced in \cite{pereira}. The latter method is a sampling-based extension of
the Nested Decomposition method \cite{birgemulti} which relies upon computing approximations of the Bellman functions.
Both techniques have been applied to a variety of real-life problems and several enhancements of these methods
have been proposed recently, see, e.g., \cite{powellbook,shapsddp} for an overview of these techniques and
\cite{philpot,lecphilgirar12,guiguessiopt2016,guiguesinexact2018,guiguesbandarra17} for convergence analysis and some of recent variations. However, to the best of our knowledge, apart from some asymptotic convergence results, no theoretical performance guaranty is available for these methods.

Another general approach which has proved successful for various classes of
uncertain optimization problems relies on
restricting the control policies to belong to certain parametric families of functions of those uncertain parameters which are known at the moment when the decision is to be applied. This approach often allows for a tractable
reformulation of the corresponding optimization problem. The simplest rules of this type are affine rules which were studied for stochastic and chance-constrained programming \cite{charnes1958cost,charnes1963deterministic} and, more recently, in the context of robust optimization
\cite{nem2003,ben2009robust}. Though affine decision strategies are generally suboptimal \cite{wetsldr} they usually provide a convenient conservativeness/numerical cost tradeoff and have been successfully used in many real-life applications, such as production management \cite{vincentrobust09},
portfolio optimization \cite{kuhmapp}, and unit commitment problems \cite{lorcaetal2015}. Nevertheless, to reduce the conservativeness of affine decisions at the price of increased computational complexity, other families of parametric decision rules were recently proposed, e.g., those using liftings as in \cite{kuhn15}, polynomial decision rules \cite{kuhm11poly},
or projected affine decision rules combined with dynamic chance constraints
as in \cite{vincentrene2017}, among others.

In this paper, we discuss a new family of decision rules for multistage uncertain problems called
{\em Constant Depth Decision Rules} (CDDRs). We consider two classes of problems: problems with discrete uncertain parameters which, at each stage, take values in a finite set, and
problems with polytopic uncertainties taking values, at each stage, in a convex hull of a finite number of points. When dealing with discrete uncertainties the corresponding decision $x_t(\cdot)$  at stage $t$ is additive with memory $\mu$, i.e., is a sum of $t$ (arbitrary) functions $u^t_s(\cdot)$, $s=1,...,t$, of $\mu$ consecutive observations of uncertain parameters preceding stage $s$  (here $\mu$ is a given {\em depth} parameter of the rule).
When uncertainties are polytopic,
decision variable  $x_t(\cdot)$ at stage $t$ is a poly-affine function, i.e. a sum of functions $v^t_s(\cdot)$, $s=1,...,t$, of $\mu$ consecutive observations of uncertainties preceding stage $s$ which are affine in each argument (see the precise definition in  Section \ref{sec:poly}). We study a class of %deterministic and stochastic
optimization problems under uncertainty with linear constraints
%with certain coefficients at the variables
and uncertain right-hand sides. Our objective is to minimize convex deterministic or stochastic objective over the set of Constant Depth Decision Rules satisfying the constraints.
Our principal contributions are the following.
\begin{itemize}
\item We provide a tractable
equivalent reformulation of the problem to be solved to compute optimal CDDRs. Specifically, for an $N$-stage optimization problem with no more than $m$ linear constraints and $n$ variables
at each stage $t$, and at most $d$ possible values of uncertain parameters at each stage $t$, and right-hand side uncertainties  of the same additive structure as  our decision rules, we rewrite the system of linear constraints on CDDRs $x_t(\cdot)$    as a linear system
of constraints
with the total numbers of scalar decision variables and scalar constraints bounded with $O(1)(m+n)N^2 d^\mu$.
\item We establish similar results for the case when the uncertainty at each stage runs through a polytope with at most $d$ vertices, and decision rules and uncertain right-hand sides of the linear constraints are allowed to be poly-affine, with memory depth $
\mu$, functions of the uncertainties.
\item Finally, to illustrate an application of the proposed approach, we compare it to that of SDDP on a toy problem of hydro-thermal production planning using a MATLAB toolbox {\tt\url{https://github.com/vguigues/Constant_Depth_Decision_Rules_Library}} for computing optimal CDDRs for multistage linear programs
with discrete uncertainties.
\end{itemize}
The paper is organized as follows. In Section \ref{discrunc} we describe the CDDR approach in the case of discrete uncertainties. In particular, we show how to reformulate the system of uncertain linear constraints of the problem of interest, when solved using CDDRs, as a linear system of inequalities on the coefficients of the CDDRs, with size of the system polynomial, the memory depth $\mu$ being fixed, in the sizes of the problem of interest. In Section \ref{sec:continuous} we develop CDDRs applied to problems with
uncertainties supported on polytopes. We illustrate the use of the proposed methodology on an example of hydro-thermal production planning application described in
Section \ref{sec:appli}. Results of numerical experiments comparing CDDRs to SDDP are reported in Section \ref{sec:num} along with a comparative discussion of these approaches to the application in question.
\section{Fixed memory decision making under discrete uncertainty}\label{discrunc}
\subsection{The problem}
Consider the situation where we control a system
$S$ evolving over stages $1,...,N$ and  affected by our decisions $x_t\in\bR^{n_t}$ and external disturbances (``uncertainties'') $\xi_t$, $1\leq t\leq N$, where $x_t$ is allowed to be a function of $\xi^t=(\xi_1,...,\xi_t)$: $x_t=x_t(\xi^t)$. We assume that\\
\begin{enumerate}
\item[{\bf A.1.}] For every $t\leq N$, disturbance  $\xi_t$ takes values in a finite set of cardinality $d_t$, which we  identify, without loss of generality,
with $D_t=\{1,2,...,d_t\}$.\\
\item[{\bf A.2.}] Feasibility of controls $x_t$ is expressed by linear constraints
\begin{equation}\label{eq00}
\sum_{\tau=1}^t A^{t\tau}x_\tau(\xi^\tau)\leq b_t(\xi^t)\in\bR^{m_t},\,t=1,...,N
\end{equation}
which should be satisfied for all trajectories $\xi^N\in D^N=D_1\times...\times D_N$.\\
\end{enumerate}

Under these restrictions, we want to minimize a given objective.
In order to make the problem computationally tractable, we further restrict the structure of the decision rules we allow for. We also consider a specific class of objective functions
and a special form of dependence of $b_t(\xi^t)$ on $\xi^t$, as explained
below.\\
\paragraph{Preliminaries: additive functions with memory $\mu$.} To avoid messy notation, we augment a sequence $\xi^\tau=(\xi_1,\xi_2,\ldots,\xi_\tau)\in D^\tau:=D_1\times...\times D_\tau$ with terms $\xi_s$ with nonpositive indices $s\leq 0$; all these terms take values in the singletons $D_s=\{1\}$, that is, from now on $\xi_s=1$ when $s\leq 0$. Besides this, for a {\sl trajectory of disturbances}---a sequence $\xi^N=(\xi_1,...,\xi_N)\in D^N$---and $p\leq q\leq N$, we denote by $\xi_{p:q}$ the fragment
$(\xi_p,\xi_{p+1},\ldots,\xi_q)$ of $\xi^N$, with the already made convention that $\xi_s=1$ for $s\leq0$. Finally, let us agree that notation $\xi^t$ when used in the same context with  $\xi_p$ or  $\xi_{p:q}$, with $q\leq t$, means that $\xi_p$ is the $p$-th entry, and $\xi_{p:q}$ is the fragment $(\xi_p,\xi_{p+1},...,\xi_q)$ of the sequence $\xi^t=(\xi_1,...,\xi_t)\in D^t$.

\begin{definition}
Given a positive integer $\mu$, let us call a function $g(\xi_1,\xi_2,...,\xi_t): D^t\to\bR^\nu$ {\sl additive with memory $\mu$}, if
$$
g(\xi_1,...,\xi_t)=\sum_{\tau=1}^t u_{\tau\xi_{\tau-\mu+1:\tau}},
$$
where {\sl coefficients} $u_{\tau\xi_{\tau-\mu+1:\tau}}$ of $g$ take values in $\bR^\nu$.
\end{definition}
\begin{example}
An additive with memory  1
 function $g:D^t\to\bR^\nu$ is specified by the collection $\{u_{\tau\xi}\in\bR^\nu:1\leq\tau\leq t,1\leq\xi\leq d_\tau\}$ of
 coefficients of $g$, and the value of $g$ at a $\xi^t\in D^t$ is the sum $\sum_{\tau=1}^tu_{\tau\xi_\tau}$ of the coefficients taken ``along the trajectory $\xi^t=(\xi_1,...,\xi_t)$.''
\end{example}

\paragraph{Structural restrictions.}
In the sequel, aside of Assumptions {\bf A.1-2} we have already made, we fix a positive integer $\mu$ and impose the following restrictions on the structure of problem constraints and decision rules:\\
\begin{enumerate}
\item[{\bf A.3.}] The right-hand sides $b_t(\xi^t)$  in design specifications (\ref{eq00}) are additive with memory $\mu$.\\
\item[{\bf A.4.}] The decision rules $x_t(\xi^t)$ are restricted to be additive with memory  $\mu$.
\end{enumerate}

According to {\bf A.3}, we have
\begin{equation}\label{wehaverhs}
b_t(\xi^t)=\sum_{\tau=1}^t\beta^t_{\tau \xi_{\tau-\mu+1:\tau}}
\end{equation}
for some collection
$$
\begin{array}{c}
\left\{\beta^t_{\tau\xi}\in\bR^{m_t},1\leq\tau\leq t, \xi  \in D_{\tau-\mu+1:\tau}\right\},\\
D_{p:q}=D_p\times D_{p+1}\times...\times D_q.\\
\end{array}
$$
Observe that A.3 encompasses a large class of right-hand sides. For instance,
\if{
\aic{$b_t$ affine in $\xi^t$ and $\xi_t$ a linear model (as in the application
considered in Section \ref{sec:appli} where $\xi_t$ is Periodic Autoregressive
with discretized noises); nonlinear models
of dependence are also allowed.}
}\fi
it allows $b_t$ which are affine in $\xi^t$ with $\xi_t$ satisfying a linear model, e.g., the autoregressive $\xi_t$'s as is the case  in the application considered in Section \ref{sec:appli}. Furthermore, $b_t(\xi^t)$ may be a nonlinear function of $\xi^t$;
for example,  for $\mu=2$ one may consider $b_t(\cdot)$ of the form
$\sum_{j=1}^t f_{t j}(\xi_j,\xi_{j-1})$
with arbitrary functions $f_{t j}(\cdot)$, etc.

Similarly, by {\bf A.4}, candidate decision rules $x^N=\{x_t(\cdot):1\leq t\leq N\}$ can be parameterized by collections
$$u^N=\left\{u^t_{\tau\xi}\in\bR^{n_t}:1\leq\tau\leq t\leq N,\xi\in D_{\tau-\mu+1:\tau}\right\}$$
according to
\begin{equation}\label{wehavexs}
x_t(\xi^t)=\sum_{\tau=1}^t u_{\tau\xi_{\tau-\mu+1:\tau}}^t.
\end{equation}
\par
Our final assumption is as follows:\\
\begin{enumerate}
\item[{\bf A.5.}] The objective to be minimized is an efficiently computable convex function $f(u^N)$ of the vector $u^N$ of parameters of a candidate decision rule.\\
\end{enumerate}
Immediate examples of the objectives of the required structure are given by the following construction: we are given a real-valued function $F(x^N,\xi^N)$ which is convex in $x^N=[x_1;x_2;...;x_N]\in\bR^{n_1}\times...\times\bR^{n_N}$, and a probability distribution $P$ on the set
$D^N$ of $N$-element trajectories of disturbances, and our objective is the expectation ${f(u^N)=}\bE_{\xi^N\sim P}\left\{F(x^N(\xi^N),\xi^N)\right\}$ of the ``loss'' $F(x^N,\xi^N)$  as evaluated at our controls.
When decision rules $x_t(\xi^t)$ are additive, with memory $\mu$, objectives of this type are convex in $u^N$; whether they are efficiently computable depends on the structure of $P$. Computability takes place when $P$ is {known and} supported on a subset of $D^N$ of moderate cardinality. {When this is not the case, but we can efficiently sample from $P$,  we} can arrive at the latter situation when {replacing} the actual objective {by} its Sample Average Approximation (that is, approximating $P$ by the uniform distribution on a reasonably large sample of trajectories of disturbances drawn from $P$, see \cite{shadenrbook}).
\par Another {important} example of objective satisfying {\bf A.5} is the maximum, over all trajectories $\xi^N\in D^N$, of a linear functional $\sum_{t=1}^N\eta_t^Tx_t(\xi^t)$ of the control trajectory, cf. {Section \ref{sec:vmod}}.
\subsection{Processing the problem}\label{sec:constr}
Treating $u^N=\{u^t_{\tau\xi}\in\bR^{n_t}:1\leq \tau\leq t\leq N, \xi\in D_{\tau-\mu+1:\tau}\}$ as our design variables,  the constraints (\ref{eq00}) read
\begin{equation}\label{emeq2}
\sum_{\tau=1}^tA^{t\tau}x_\tau(\xi^\tau)-b_t(\xi^t)\leq 0\;\forall (t,1\leq t\leq N,\xi^N\in D^N),
\end{equation}
or equivalently (see (\ref{wehaverhs}), (\ref{wehavexs}))
\begin{equation}\label{emeq3}
\sum\limits_{s=1}^t\left[\sum\limits_{\tau=s}^t A^{t\tau} u^\tau_{s\xi_{s-\mu+1:s}}-\beta^t_{s\xi_{s-\mu+1:s}}\right]\leq 0\;\forall (t,1\leq t\leq N,\xi^N\in D^N).
\end{equation}
The crucial fact for us is that {\sl constraints {\rm (\ref{emeq3})} can be {reduced to} an explicit system of linear inequality constraints on the design variables $u^t_{\tau\xi}$ and additional ``analysis'' variables}. The construction goes as follows.\\
\par
{\bf 1.} We introduce variables $y^t_{\tau\xi}\in\bR^{m_t}$, $1\leq \tau\leq t\leq N, \xi\in D_{\tau-\mu+1:\tau}$, and link them to our decision variables $u^t_{\tau\xi}$ by linear equality constraints
\begin{equation}\label{emeq100}
y^t_{s\xi}=\sum\limits_{\tau=s}^t A^{t\tau} u^\tau_{s\xi}-\beta^t_{s\xi},\,
\forall \xi\in D_{s-\mu+1:s},1\leq s\leq t\leq N.
\end{equation}
In terms of these variables (\ref{emeq3}) reads
\begin{equation}\label{emeq90}
\sum\limits_{s=1}^ty^t_{s\xi_{s-\mu+1:s}}\leq0\,\,\,\forall \xi^N\in D^N.
\end{equation}
To avoid messy notation, we describe our subsequent actions separately for the case of $\mu=1$ and of $\mu>1$.\\
\par
{\bf 2.A: {Case of} $\mu=1$.} For every $t\in\{1,...,N\}$, we introduce variables $z^t_s\in\bR^{m_t}$, $1\leq s\leq t$, and impose the linear inequalities
%\begin{equation}\label{emeq101a}
%\begin{array}{rcll}
%z^t_s&\geq& y^t_{s\xi}+z^t_{s+1}\forall \xi\in D_s,s=t,t-1,...,1,&(a)\\
%z^t_1&\leq& 0&(b)
%\end{array}
%\end{equation}
\begin{subequations}
\begin{align}
z^t_s&\geq y^t_{s\xi}+z^t_{s+1}\forall \xi\in D_s,s=t,t-1,...,1,\label{emeq101aa}\\
z^t_1&\leq 0\label{emeq101ab}
\end{align}
\label{emeq101a}
\end{subequations}
where $z^t_{t+1}\equiv 0$.
Clearly, {the} $i$-th entry $[z^t_s]_i$ in $z^t_s$ is an upper bound on
$\max_{ \xi^N \in D^N   }\left[\sum_{r=s}^t y^t_{r\xi_{r}}  \right]_i$, and constraints
\eqref{emeq101aa} allow to make this bound equal to the latter quantity. Consequently, the system $\cS$ of  constraints (\ref{emeq100}){ and}
(\ref{emeq101a}) on variables $u$, $y$, $z$ provides a polyhedral representation of the solution set of (\ref{emeq2}). In other words,  a collection $u^N$ of  actual design variables $u^t_{\tau\xi}$ satisfies constraints (\ref{emeq2}) if and only if $u^N$ can be extended,
by properly selected values of $y$- and $z$-variables, to a feasible solution of $\cS$.
On the other hand, (\ref{emeq101a}) is entrywise decomposable (it is a collection of $m_t$ systems of linear inequalities, with {the} $i$-th system involving only {the} $i$-th entries in $y$-and $z$-vectors), and as far as the $i$-th entries in $y^t_{s\xi}$ and $z^t_s$ are concerned, (\ref{emeq101aa}) is nothing but the ``backward'' Dynamic Programming description of an upper bound $[z^t_1]_i$ on $\max\limits_{\xi^N\in D^N}\left[\sum_{s=1}^t y^t_{s\xi_{s}}\right]_i$ (recall that we are in the case of $\mu=1$), while (\ref{emeq101ab}) says that the resulting bound should be nonpositive for all $i$, exactly as required in (\ref{emeq90}).\\
\par{\bf 2.B: {Case of} $\mu\geq 2$.} As in the   case $\mu=1$, what follows is nothing but backward Dynamic Programming description, expressed by linear inequalities,  of vector
$z^t_{1\xi_{2-\mu:0}}$ with $i$-th entry, $i\leq m_t$, upper-bounding $\max\limits_{\xi^N\in D^N}\left[\sum_{s=1}^t y^t_{s\xi_{s-\mu+1:s}}\right]_i$. As is immediately seen, to get this description it suffices to introduce variables $z^t_{s\eta}\in\bR^{m_t}$, $1\leq s\leq t$, $\eta\in D_{s-\mu+1:s-1}$, and subject them, along with the $y$-variables, to linear constraints
\begin{equation}\label{emeq101b}
\begin{array}{lcl}%{lcll}
z^t_{t \xi_{t-\mu+1:t-1}} &\geq& y^t_{t\xi_{t-\mu+1:t}}\forall\xi_t \in D_t,\\
z^t_{s\xi_{s-\mu+1:s-1}}&\geq&y^t_{s\xi_{s-\mu+1:s}}+z^t_{(s+1)\xi_{s-\mu+2:s}}\,\,\forall \xi_{s-\mu+1:s} \in D_{s-\mu+1:s},s=t-1,t-2,...,1,\\
z^t_{1\xi_{2-\mu:0}}&\leq&0.
\end{array}
\end{equation}
Similarly to the case of $\mu=1$, the system of all constraints (\ref{emeq100}) and (\ref{emeq101b}) on variables $u$, $y$, $z$ gives a polyhedral representation of the solution set of (\ref{emeq2}).
 \par
The bottom line is that {\sl under Assumptions {\bf A.1-5}, the problem of interest can be straightforwardly reduced to the problem of minimizing an efficiently computable convex objective $f(u^N)$ over $u$-, $y$-, $z$-variables satisfying an explicit system of linear inequality constraints.} Note that for every fixed $\mu$, the total number of variables and constraints in the resulting problem $\cP$ is polynomial in the sizes of the problem of interest. Specifically, assuming $d_t\leq d$, $m_t\leq m$ and $n_t\leq n$ for all $t$, the total numbers of scalar decision variables and scalar constraints in $\cP$ do not exceed $O(1)(m+n)N^2d^\mu$.
\par Finally, observe that
if $b_t$ is additive with memory
$\mu_t$ then it is also additive with memory
$\mu$ for any $\mu \geq \max_{1 \leq i \leq N} \mu_i$ and we can
apply the proposed methodology  utilizing additive with memory
$\mu$ decision rules for any $\mu \geq \max_{1 \leq i \leq N} \mu_i$.

\subsection{Modifications}\label{sec:vmod}
In the above exposition, we treated  vectors $\beta^t_{\tau\xi}$ as part of the data. It is immediately seen that when changing the status of some of $\beta$'s from being part of the  data to being additional decision variables and adding, say, linear constraints on these variables, we preserve tractability of the resulting problem: our backward Dynamic Programming still allows us to convert constraints (\ref{emeq100}) and (\ref{emeq90}) modified in this way into an explicit system of linear inequalities on ``variables of interest'' (components of $u^N$ and new design variables coming from $\beta$'s) and additional $y$- and $z$-analysis variables.

An immediate application of this observation is as follows. Suppose that Assumptions {\bf A.1-4} hold and assume that our goal is to minimize the worst case (i.e., the largest over $\xi^N\in D^N$) value of the function
$$
F[\{x_t(\xi^t)\}_{t=1}^N]=\max_{\ell\leq L}\sum_{t=1}^Nh_{t\ell}^Tx_t(\xi^t).
$$
To this end it suffices to augment our original design variables $u^N=\{u^t_{\tau\xi}\}$ parameterizing additive, with memory $\mu$, candidate decision rules  with a new decision variable $w$ and extend the original system
$\sum_{\tau=1}^NA^{N\tau}x_\tau(\xi^\tau)\leq b_N(\xi^N)$ of the last stage constraints by adding to it constraints
$$
\sum_{t=1}^Nh_{t\ell}^Tx_t(\xi^t)\leq w,\,\ell=1,...,L.
$$
As a result, we get a parametric, the parameter being $\xi^N\in D^N$, system of linear inequalities on $u^N$ and $w$. Applying backward Dynamic Programming in exactly the same way as above, we convert this system into an explicit system of linear inequalities on $u^N$, $w$, and additional $y$- and $z$-variables. Minimizing the worst-case value of the above criterion is thus reduced to an explicit Linear Programming problem in $w$ and $u$-, $y$-, $z$-variables.
\par Until now we have assumed that control feasibility is expressed in terms of the system \rf{eq00} of linear constraints with uncertain right-hand sides. It may be worth mentioning that the proposed approach can be straightforwardly modified to deal with linear constraints with uncertain matrices or even specific nonlinear constraints at the price of restricting severely  the class of control strategies.
\par
Indeed, let us assume from now on that the control action $x_t$ depends solely on  $\xi_{t-\mu+1:t}=(\xi_{t-\mu+1},..., \xi_{t})$ so that representation \rf{wehavexs} reduces to
\begin{equation}\label{simpledecisionrules}
x_t(\xi^t)=u^t_{\xi_{t-\mu+1:t}}.
\end{equation}
Now, let us consider an uncertain linear system \rf{eq00} satisfying Assumptions {\textbf{A.1.}}, {\textbf{A.2.}} and {\textbf{A.3.}}
with uncertain technology matrices $A^{t \tau}$. More precisely, we assume that matrices
$A^{t\tau}$ depend on the fragment
$\xi_{\tau-\mu+1:\tau}=(\xi_{\tau-\mu+1},\xi_{\tau-\mu+2},\ldots,\xi_\tau)$ of $\xi^\tau$, i.e., constraints \eqref{eq00} are replaced with the constraints
\[
\forall (t:1\leq t\leq N):\; \sum_{\tau=1}^t A^{t\tau}(\xi_{\tau-\mu+1:\tau}) u^\tau_{\xi_{\tau-\mu+1:\tau}} \leq b_t(\xi^t),
\]
where $b_t$, as before, is additive with memory $\mu$.

As was done in Section \ref{sec:constr}, we can rewrite the constraints replacing $y$-variables with
$$
y^{t}_{s \xi}=\displaystyle
A^{t s}(\xi) u^s_{\xi}-\beta^{t}_{s \xi},\;\;
     \xi \in D_{s-\mu+1:s},1\leq s\leq t \leq N,
$$
and therefore all the machinery developed in Section \ref{sec:constr} can be applied.\\
Given $x_0\in \bR^{n_0}$, let us now consider a system of convex nonlinear constraints
\begin{equation}\label{constraintsnonconvex}
\begin{array}{l}
\displaystyle \sum_{\tau=1}^{t} G^{t \tau}(x_{\tau}, x_{\tau-1},
\xi_{\tau-\mu+1:\tau}) \leq 0 \in\bR^{m_t},\,t=1,...,N,
\end{array}
\end{equation}
``coupling'' control actions at subsequent stages, which should be satisfied for all trajectories $\xi^N\in D^N=D_1\times ... \times D_N$
with $(\xi_t)$ satisfying  {\textbf{A.1}.}
We suppose that components
$G_i^{t \tau}(x,x',\xi_{\tau-\mu+1:\tau}), \,i=1,\ldots,m_t,$ of $G^{t \tau}$
are convex in $x$ and $x'$  for all possible $(\xi_1,\xi_2,\ldots,\xi_N)$ and all $t,\tau$.
As above, using decision rules \eqref{simpledecisionrules}, we derive the following
representation of constraints \eqref{constraintsnonconvex} in variables $u^t_{\xi_{t-\mu+1:t}}$, $t=1,\ldots,N, \;\xi_{t-\mu+1:t} \in
D_{t-\mu+1:t}$, and $z_{s,\xi_{s-\mu:s-1}}^{t}$, $2 \leq s \leq t \leq N$, $\xi_{s-\mu:s-1} \in D_{s-\mu:s-1}$, by Dynamic Programming:
\bse
&&G^{1 1}(u^1_{\xi_{2-\mu:1}},x_0,\xi_{2-\mu:1})\leq 0,\\
&&\forall t,\, 2\leq t\leq N, \;\forall \xi_{t-\mu:t} \in D_{t-\mu:t}:\\
&&\displaystyle
z_{t,\xi_{t-\mu:t-1}}^{t} \geq G^{t t}
(u^t_{\xi_{t-\mu+1:t}},u^{t-1}_{\xi_{t-\mu:t-1}},\xi_{t-\mu+1:t}),\\
&&\forall s,t,\, 3 \leq s \leq t \leq N,\;\forall \xi_{s-\mu-1:s-1} \in D_{s-\mu-1:s-1}:\\
&&z_{s-1,\xi_{s-1-\mu:s-2}}^{t} \geq G^{t\, s-1}
(u^{s-1}_{\xi_{s-\mu:s-1}},u^{s-2}_{\xi_{s-\mu-1:s-2}},
\xi_{s-\mu:s-1})+ z_{s,\xi_{s-\mu:s-1}}^{t},\\
&&\forall t,\, 2\leq t\leq N:\\
&&z_{2,\xi_{2-\mu:1}}^{t} + G^{t 1}( u^1_{\xi_{2-\mu:1}},x_0,\xi_{2-\mu:1}) \leq 0.
\ese
%with given $x_0$ and
%$\xi_{s-\mu:s-1} \in D_{s-\mu:s-1}$.
%
%
\section{Fixed memory decision making under polytopic uncertainty}\label{sec:continuous}
\subsection{The problem}\label{sec:poly} So far, we have considered multi-stage decision making under discrete uncertainty, where the external disturbance acting at the controlled system at time $t$ takes one of $d_t$ values known in advance. Let us now consider the case of {\sl polytopic} uncertainty, where the disturbance at time $t$ is a vector $\zeta_t$ taking values in a given polytope
$\Delta_t \subset \bR^{\nu_t-1}$. As before,  we allow  for our decision at time $t$, $x_t\in\bR^{n_t}$, to depend on the sequence $\zeta^t=(\zeta_1,...,\zeta_t)$. From now on we make the following assumptions (cf. Assumptions {\bf A.1-2}):\\
 \begin{enumerate}
\item[{\bf B.1.}] For every $t\leq N$, disturbance  $\zeta_t$ takes values in polytope $\Delta_t\subset\bR^{\nu_t-1}$ given by the list of  $d_t$ scenarios $\chi_{ts}$:
\begin{equation}\label{deltat}
\Delta_t=\Conv\{\chi_{ts},1\leq s \leq d_t\}\subset\bR^{\nu_t-1}.
\end{equation}

We also assume that the scenarios affinely span $\bR^{\nu_t-1}$ (this is w.l.o.g., since we can always replace the embedding space $\bR^{\nu_t-1}$  of $\Delta_t$ by the affine span of $\Delta_t$). Thus, $\nu_t\leq d_t$ for all $t$.
\item[{\bf B.2.}] Feasibility of controls $x_t$ is expressed by linear constraints
\begin{equation}\label{eq00n}
\sum_{\tau=1}^t A^{t\tau}x_\tau(\zeta^\tau)\leq b_t(\zeta^t)\in\bR^{m_t},\,t=1,...,N
\end{equation}
which should be satisfied for all trajectories $\zeta^N\in \Delta^N=\Delta_1\times...\times \Delta_N$.\\
\end{enumerate}
\par Under these restrictions, we want to minimize a given objective. To ensure computational tractability, as in the case of discrete disturbances,  we impose structural restrictions on the allowed $x_t(\cdot)$'s, $b_t(\cdot)$'s, and on the objective. Our main restriction is that the policies $x_t(\cdot)$ and the right-hand sides $b_t(\cdot)$ are {\sl poly-affine with memory $\mu$}.\\
\paragraph{Poly-affine functions with memory $\mu$.} For notational convenience, we augment a sequence $\zeta^\tau=(\zeta_1,\zeta_2,...,\zeta_\tau)\in\Delta^\tau:=\Delta_1\times...\times\Delta_\tau$ with terms $\zeta_s$ with nonpositive indices $s\leq 0$; all these terms take values in the singletons $\Delta_s=\{0\}=\bR^0$, that is, from now on $\zeta_s=0\in\bR$ when $s\leq 0$.

\begin{definition}
Given a positive integer $\mu$, let us call function $g(\zeta_1,\zeta_2,...,\zeta_t): \Delta^t\to\bR^\nu$ {\sl poly-affine with memory $\mu$}, if
$$
g(\zeta_1,...,\zeta_t)=\sum_{\tau=1}^t g_\tau(\zeta_{\tau-\mu+1},...,\zeta_\tau),
$$
where {{\sl every component} $g_\tau(\zeta_{\tau-\mu+1},...,\zeta_\tau)$ of $g$ takes values in $\bR^\nu$ and is affine in each of its arguments $\zeta_{\tau-\mu+1},...,\zeta_\tau$}.
\end{definition}
\begin{example}
\begin{itemize}
\item An affine vector-valued function $g(\zeta^t)=\sum_{\tau=1}^t  [p_\tau +P_\tau\zeta_\tau]$ of $\zeta_1,...,\zeta_t$ is poly-affine with memory $\mu=1$; its components are $g_\tau(\zeta_\tau)=p_\tau+P_\tau\zeta_\tau$.
\item A general poly-affine vector-valued function $g(\zeta^t)$ of $\zeta_\tau\in\bR^{\nu_\tau}$, $1\leq\tau\leq t$, with memory $\mu=2$ is a function representable in the form
$$
g(\zeta^t)=\sum_{\tau=1}^tg_\tau(\zeta_{\tau-1},\zeta_\tau),\,\,g_\tau(\zeta_{\tau-1},\zeta_\tau)=p_\tau+P_\tau\zeta_\tau+\sum\limits_{{1\leq r\leq\nu_{\tau-1},\atop1\leq s\leq \nu_\tau}}[\zeta_{\tau-1}]_r[\zeta_{\tau}]_s
f_{\tau,r,s},
$$
where $p_\tau,f_{\tau,r,s},P_\tau$ are vectors and matrix of appropriate sizes, and $[a]_i$ stands for $i$-th entry of a vector.
\end{itemize}
\end{example}

Given $t\leq N$, we denote by $\lambda(\zeta_t)\in\bR^{\nu_t}$ the
coordinates of $\zeta_t \in \bR^{\nu_t-1}$
in an affine basis of $\bR^{\nu_t-1}$. For instance,
if the affine basis is made of
the standard basis vectors and the origin, we get
$$
\lambda_i(\zeta_t)=[\zeta_t]_i,1\leq i<\nu_t,\,\lambda_{\nu_t}(\zeta_t)=1-\sum_{i=1}^{\nu_{t}-1}[\zeta_t]_i.
$$
We also set $\lambda_1(\zeta_t)\equiv 1$ when $t\leq 0$ (and, according to our convention,
$\zeta_t=0\in\Delta_0 =\acute{}\bR^0 =\{0\}$). From affinity of $g_\tau(\zeta_{\tau-\mu+1},...,\zeta_\tau)$ in each argument it follows that
$$
\begin{array}{c}
g_\tau(\zeta_{\tau-\mu+1},...,\zeta_\tau)=\sum_{\varkappa\in \cI_\tau}\left[\prod\limits_{s=1}^\mu\lambda_{\varkappa_s}(\zeta_{\tau-\mu+s})\right]g_{\tau\varkappa},\\
\cI_\tau=\{\varkappa=(\varkappa_1,...,\varkappa_\mu): 1\leq \varkappa_s\leq \nu_{\tau-\mu+s},1\leq s\leq \mu\}.\\
\end{array}
$$

As a result, a poly-affine function of $\zeta^t$ taking values in $\bR^\nu$ is fully specified by the {\sl collection of its coefficients}
$$
g^t=\left\{g_{\tau\varkappa}\in\bR^\nu:1\leq \tau\leq t,\varkappa \in \cI_\tau\right\}
$$
according to
\begin{equation}\label{multiaffine}
\begin{array}{c}
g(\zeta^t)=\sum_{\tau=1}^t\sum_{\varkappa\in\cI_\tau}\left[\prod\limits_{s=1}^\mu\lambda_{\varkappa_s}(\zeta_{\tau-\mu+s})\right]g_{\tau\varkappa}.
\end{array}
\end{equation}
Note that every collection $\{g_{\tau\varkappa}\in\bR^\nu:1\leq\tau\leq t,\varkappa\in\cI_\tau\}$ is a collection of coefficients of a poly-affine with memory $\mu$ function $g(\zeta^t)$ taking values in $\bR^\nu$.\\
\paragraph{Structural restrictions} imposed in the sequel on the decision rules and right-hand sides in the constraints (\ref{eq00n}) are as follows (cf. {\bf A.3-4}):\\
\begin{enumerate}
\item[{\bf B.3.}] The right-hand sides $b_t(\zeta^t)$  in design specifications (\ref{eq00n}) are poly-affine with memory $\mu$.\\
\item[{\bf B.4.}] The decision rules $x_t(\zeta^t)$ are restricted to be poly-affine with memory $\mu$.\\
\end{enumerate}
By {\bf B.4}, candidate decision rules $x^N=\{x_t(\cdot):1\leq t\leq N\}$ in question can be parameterized by finite-dimensional collections
\begin{equation}\label{eqvN}
v^N=\{v^t_{\tau\varkappa}\in\bR^{n_t}:1\leq \tau\leq t\leq N,\varkappa\in\cI_\tau\}
\end{equation}
according to
\begin{equation}\label{xtzeta}
x_t(\zeta^t)=\sum_{\tau=1}^t \sum_{\varkappa\in\cI_\tau}\left[\prod\limits_{s=1}^\mu \lambda_{\varkappa_s}(\zeta_{\tau-\mu+s})\right]v^t_{\tau\varkappa}.
\end{equation}
For any selection of vectors $v^t_{\tau\varkappa}\in\bR^{n_t}$ in (\ref{eqvN}), the resulting collection specifies candidate decision rules $x_t(\cdot)$, $t\leq N$, satisfying {\bf B.4}.\\
\par Finally, we make the following assumption (cf. {\bf A.5}):\\
\begin{enumerate}
\item[{\bf B.5.}] The objective to be minimized is an efficiently computable convex function $f(v^N)$ of the vector $v^N$ of parameters of a candidate decision rule.\\
\end{enumerate}
\subsection{Processing the problem}\label{sec:constrcont}
Let $D_t=\{1,...,d_t\}$ be the set of indices of scenarios $\chi_{ts}$ specifying $\Delta_t$ according to (\ref{deltat}). Our objective in this section is to reduce the ``continuous'' problem posed in Section \ref{sec:poly}
to a linear optimization problem (with, as for the discrete case,
a total number of scalar decision variables and scalar constraints not exceeding $O(1)(m+n)N^2d^\mu$ when $m_t \leq m$, $n_t \leq n$, $d_t \leq d$). To this end, we show
that poly-affine decisions  $x_t(\zeta^t)$ of the form (\ref{xtzeta}) satisfy the constraints \rf{eq00n} on every trajectory
$\zeta^N\in\Delta^N$ of the ``continuous'' uncertainty if and only if they satisfy the constraints for all
{\sl scenario} trajectories $\zeta^t$, those   which are sequences of {scenarios} $\chi_{t\xi_t}$, $\xi_t\in D_t,1\leq t\leq N$.
As we shall see, this simple observation allows us to reduce the continuous problem to slightly modified discrete problem from Section \ref{discrunc}. An informal outline of the reduction is as follows: there is a natural way to restrict the continuous problem onto the scenario trajectories of uncertainty, thus arriving at a discrete problem from Section \ref{discrunc}, with uncertainties $\xi_t\in D_t$ stemming from  the indices of scenarios $\chi_{t\xi_t}$. With this reduction, the restrictions of poly-affine, with memory depth $\mu$, control policies $\{x_t(\zeta^t)\}_{t\leq N}$ of continuous problem onto the scenario trajectories are exactly policies of the form (\ref{wehavexs}) for the discrete problem, with vectors  $u_{\tau\xi_{\tau-\mu+1:\tau}}^t$ being linear images, under known linear mappings, of the collections of coefficients of  poly-affine functions $\{x_t(\xi^t\}_{t\leq N}$. Because, as we have mentioned, a poly-affine control policy for the continuous problem is feasible if and only if its restriction on the scenario  trajectories is feasible for the discrete problem, these observations reduce the continuous problem to the discrete one.
\par
To construct the discrete problem it is convenient to associate
to every ``trajectory of indices'' $\xi^t\in D^t:=D_1\times...\times D_t$,
a trajectory of disturbances
$$
\zeta^t[\xi^t]=\{\zeta_\tau[\xi_\tau]:=\chi_{\tau\xi_\tau}:1\leq \tau\leq t\}\in\Delta^t:=\Delta_1\times...\times \Delta_t.
$$
Clearly, when $\xi^\tau\in D^\tau$ is the initial fragment of $\xi^t\in D^t$, then $\zeta^\tau[\xi^\tau]$ is the initial fragment of $\zeta^t[\xi^t]$.  Let us make two immediate observations:
\begin{proposition}\label{lem1} Let $g(\zeta^t)$ be a poly-affine, with memory $\mu$, function taking values in $\bR^\nu$:
\begin{equation}\label{gis}
\begin{array}{c}
g(\zeta^t)=\sum_{\tau=1}^t\sum_{\varkappa\in\cI_\tau}\left[\prod\limits_{s=1}^\mu\lambda_{\varkappa_s}(\zeta_{\tau-\mu+s})\right]g_{\tau\varkappa}.
\end{array}
\end{equation}
Then the mapping
$$
\overline{g}(\xi^t):=g(\zeta^t[\xi^t]): D^t\to\bR^\nu
$$
is additive with memory $\mu$.
\end{proposition}
{\bf Proof.}  Setting, as in Section \ref{discrunc}, $D_{p:q}=D_p\times D_{p+1}\times...\times D_q$ and given $t,\tau,\xi$ with $1\leq\tau\leq t\leq N$ and
$\xi=(\xi_{\tau-\mu+1},\xi_{\tau-\mu+2},...,\xi_\tau)\in D_{\tau-\mu+1:\tau}$, let us put
$$
{\bar{g}_{\tau\xi}}=\sum_{\varkappa\in\cI_\tau}\left[\prod\limits_{s=1}^\mu\lambda_{\varkappa_s}(\chi_{\tau-\mu+s,\xi_{\tau-\mu+s}})\right]g_{\tau\varkappa}
$$
(the right-hand side indeed depends only on $\tau$ and $(\xi_{\tau-\mu+1},\xi_{\tau-\mu+2},...,\xi_\tau)\in D_{\tau-\mu+_1:\tau}$). It remains to note that by (\ref{gis}) we have
$$
\overline{g}(\xi^t)=g(\chi_{1\xi_1},\chi_{2\xi_2},...,\chi_{t\xi_t})=\sum_{\tau=1}^t{\bar{g}_{\tau\xi_{\tau-\mu+1:\tau}}}\,\,\forall\xi^t\in D^t=D_1\times...\times D_t,
$$
as it should be for an additive with memory $\mu$ function.
\qed

\par
Our second observation is as follows:
\begin{proposition}\label{lem2} Let $f(\zeta^N)$ be a real-valued affine in every component $\zeta_t$ of $\zeta^N$
function (e.g., a poly-affine real-valued function with memory $\mu$). Then among maximizers of $f(\zeta^N)$ over $\zeta^N\in\Delta^N$ there are those of the form $\zeta[\xi^N]$ with properly selected $\xi^N\in D^N$.
\end{proposition}
{\bf Proof.}  Let $\bar{\zeta}^N$ be a maximizer of $f(\zeta^N)$ on $\Delta^N$ with the largest possible number, {let it be} $M$, of scenario  components $\bar{\zeta}_t$ (i.e., those belonging to
$\{\chi_{ts},1\leq s\leq d_t\}$). All we need to prove is that $M=N$. This is evident: assuming
 that $M<N$, i.e., that  for some $t$ the component $\bar{\zeta}_t$ of $\bar{\zeta}^N$ is not in the set $\{\chi_{ts},1\leq s\leq d_t\}$, let us ``freeze''  in $f(\zeta^N)$ all arguments $\zeta_s$ with $s\neq t$  at the values  $\bar{\zeta}_s$ and vary the $t$-th argument. Since $f$
is affine in every  $\zeta_s$, among the maximizers of the resulting function of $\zeta_t$ over $\zeta_t\in\Delta_t$ there will be an extreme point of $\Delta_t$, that is, a point from the set of scenarios of stage $t$. Replacing in $\bar{\zeta}^N$ the component $\bar{\zeta}_t$ with this scenario, we get another maximizer of $f$ on $\Delta^N$ with more than $M$ scenario components, which is impossible.
\qed

\par Now we are ready to explain how to process the problem of  interest numerically. Let us associate with our problem  (call it continuous) a {\sl discrete} problem
as follows.  The structure of the discrete problem is as considered in Section \ref{discrunc}, with inherited from the continuous problem number of stages $N$, matrices $A^{t\tau}$, and cardinalities $d_t$ of the sets $D_t$ of values of disturbance $\xi_t$ at  stage $t$. As about the right-hand sides $\overline{b}_t(\xi^t)$ in the constraints of the discrete problem, we specify them as
$$
\overline{b}_t(\xi^t)=b_t(\zeta^t[\xi^t]).
$$
Now, candidate decision rules $x_t(\zeta^t)$, $t\leq N$, in the continuous problem induce candidate decision rules
$$
\overline{x}_t(\xi^t):=x_t(\zeta^t[\xi^t])
$$
in the discrete problem. By
Proposition \ref{lem1}, restrictions {\bf B.3-4} on the structure of  $b_t(\cdot)$'s and $x_t(\cdot)$'s  ensure the validity of {\bf A.3-4} for $\overline{b}_t(\cdot)$'s and $\overline{x}_t(\cdot)$'s. Besides this, Proposition \ref{lem2} says that under restrictions {\bf B.3-4} decision rules $x_t(\cdot)$ are feasible for the continuous problem if and only if the decision rules $\overline{x}_t(\cdot)$ are feasible for the discrete problem. As we remember, the latter is equivalent to the fact that the collection
$$u^N=\{u^t_{\tau\xi}:1\leq \tau\leq t\leq N,\xi\in D_{\tau-\mu+1:\tau}\}$$
of parameters of the additive, with memory $\mu$, decision rules $\overline{x}^N$ can be augmented by properly selected $y$- and $z$-variables to yield a feasible solution to certain  system $\cS$ of linear constraints. For $\mu$ fixed, the number of constraints and variables in $\cS$, as well as the computational effort to build this system,  is polynomial in all sizes of the problem (for details, see Section \ref{discrunc}). We are in the situation where $u^N$ is obtained  from the ``primitive'' design variables, specifically, the collection $v^N$ of parameters specifying the decision rules $x_t(\cdot)$ by known to us linear transformation:
\begin{equation}\label{v2u}
\begin{array}{l}
u^t_{\tau\xi}=\sum_{\varkappa\in\cI_\tau}\left[\prod\limits_{s=1}^\mu\lambda_{\varkappa_s}(  \chi_{\tau-\mu+s,\xi_{\tau-\mu+s}})\right]v^t_{\tau\varkappa},\\
\mbox{for }1\leq\tau\leq t,\xi=(\xi_{\tau-\mu+1},\xi_{\tau-\mu+2},...,\xi_\tau)\in D_{\tau-\mu+1:\tau}.
\end{array}
\end{equation}
Extending $\cS$ to the system of linear constraints by adding variables $v^N$ and constraints (\ref{v2u}) linking the $v$- and the $u$-variables, we get a system $\cS^+$
of linear constraints in ``actual'' design variables $v^N$ and additional analysis variables (specifically, $u^N$ and  $y$- and $z$-variables inherited from $\cS$).
The bottom line is that  under Assumption {\bf B.5}, {\sl for $\mu$ fixed}, the problem of interest can be reduced to the problem of minimizing an efficiently computable convex function of
$v$-variables under a system $\cS^+$ of linear constraints on $v,u,y,z$-variables, with the total number of variables and constraints in $\cS^+$ and the computational effort of building this system which are polynomial in the sizes of the problem of interest.
\par
{Note that, similarly to the case of discrete uncertainty, an objective which is the worst-case value over all trajectories $\zeta^N$ of a linear functional $\sum_{t=1}^N\eta_t^Tx_t(\xi^t)$ of the control trajectory satisfies {\bf B.5}, cf. the beginning of  Section \ref{sec:vmod}.
Furthermore, in the case of random disturbances, an objective of the form
\[
f(v^N)=\bE_{\zeta^N\sim P}\left\{F(x^N(\zeta^N),\zeta^N)\right\}
\]
with
function $F(x^N,\zeta^N)$ which is
convex
in decision variables  $x^N=[x_1;x_2;...;x_N]$
and efficiently computable, can be replaced by its Sample Average Approximation, see, e.g. \cite{shadenrbook}.
\[
\widehat f_N(v^N)={1\over L}\sum_{\ell=1}^L F(x^N(\zeta^{N\ell}),\zeta^{N\ell})\]
over a large number $L$ of scenarios---realizations of disturbance trajectories
$\zeta^{N \ell}=(\zeta_1^{\ell},\zeta_2^{\ell},\ldots,\zeta_N^{\ell}), \ell=1,\ldots,L$.
In this case, Assumption {\bf B.5} holds for the approximate objective $\widehat f_N(v^N)$
with decision $x_t(\zeta^{t \ell})$ at stage $t$ of scenario $\ell$ linked to $v^N$-variables by the explicit
linear relation
$$
x_t(\zeta^{t \ell})=
\sum_{\tau=1}^t \sum_{\varkappa\in\cI_\tau}\left[\prod\limits_{s=1}^\mu \lambda_{\varkappa_s}(\zeta_{\tau-\mu+s}^{\ell})\right]v^t_{\tau\varkappa},
$$
see \eqref{xtzeta}.
}
\section{An application to hydro-thermal production planning}\label{sec:appli}
In this section, we illustrate the application of our methodology on a toy example of a hydro-thermal production planning problem formulated as
a Multistage Stochastic Linear Program with  linear constraints.
\subsection{Problem description}\label{sec:appli1}
Our problem modeling is as follows. Consider a set of thermal electricity production plants and hydroelectric plants distributed in $K$
regions which have to produce electricity to satisfy the
demand in each region
and each stage $t=1,\ldots,N$ of a given planning horizon. We will assume that in each region,
all thermal facilities are aggregated into a single thermal plant, and similarly,
all hydroelectric plants and reservoirs are aggregated into a single hydroelectric plant and a single reservoir.
The objective is to minimize the expected  production cost which is a sum of the cost of thermal generation and the penalties paid for the unsatisfied demand
over the planning horizon under constraints of demand satisfaction, minimal and maximal levels of the hydroelectric reservoirs and capacity constraints
of the production units.
\par
We use the following notation for time $t=1,\ldots,N$:
\begin{itemize}
\item $v_t\in \bR^K$ for reservoir levels at the end of stage $t$;
\item $w_t\in \bR^K$ for thermal generation at stage $t$;
\item $h_t\in \bR^K$ for hydroelectric generation;
\item $r_t\in \bR^K$ for unsatisfied demand;
\item $\I_t\in \bR^K$ for inflows;
\item $G_t\in \bR^{K\times K}$ is a diagonal matrix; $G_t \I_t$ is the vector of actual inflows to the reservoirs
and $(I - G_t)\I_t$ is the part of inflows automatically converted into energy by run-of-river plants;
\item $\delta_t\in \bR^K$ is the deterministic vector of energy demands;
\item $c_t\in \bR^K$ is the vector of thermal generation unit costs, and $p_t\in \bR^K$ is the vector of penalties  for the
unsatisfied demand at time $t$.
\end{itemize}
In this notation, water inventories, releases, and inflows are expressed in energy units. The hydrothermal production planning problem consists in minimizing the expected cost
\[
\bE_{\I^N}\Big[ \sum_{t=1}^N c_t^Tw_t+\sum_{t=1}^N p_t^T r_t\Big],
\]
under the following system of constraints to be satisfied almost surely:
\be
\begin{array}{rcll}
h_t&\leq& v_{t-1}-v_t+G_t \I_t,&\hbox{\mbox{[water balance]}}\\
h_t+w_t+r_t&\geq& \delta_t -(I-G_t)\I_t,&\hbox{\mbox{[demand satisfaction]}}\\
\underline{v}_t&\leq& v_t\,\leq\, \overline{v}_t,&\hbox{\mbox{[bounds on reservoir levels]}}\\
0&\leq &h_t\,\leq\, \overline{h}_t,&\hbox{\mbox{[hydroelectric generation capacity]}}\\
0&\leq &w_t\,\leq\, \overline{w}_t,&\hbox{\mbox{[thermal generation capacity]}}\\
r_t&\geq&0.&\hbox{\mbox{[nonnegativity of unsatisfied demand]}}
\end{array}
\ee{htsystem}
We assume that inflows $\I_t$, $t\leq 1$ are deterministic (i.e., $...,\I_0,\,\I_1$ are known at $t=1$ when production plan is computed for stages $t=1,\ldots,N$) and for $t \geq 2$ inflows satisfy the periodic autoregressive model:
\be
\begin{array}{rcl}\I_t&=&\theta_t+\eta_t,\\
\eta_t&=&\sum_{j=1}^{\ell_t} B_t^j \eta_{t-j}+ C_t\zeta_t
\end{array}
\ee{parmodel}
where $\theta_t\in \bR^K$ are given along with ${K\times K}$  matrices $B_t^j,\,C_t$,
while disturbances  $\zeta_t,\,,t=2,...,N,$ are independent with known distributions $P_t$ supported on the sets
$\{\chi_{t1},\ldots,\chi_{t d_t}\}$.
\par Note that inflows $\I_t$ satisfying recursive equations \eqref{parmodel}  can be straightforwardly rewritten
in the form
\be
\I_t = \nu_t + \sum_{s=2}^t R_s^t \zeta_s,
\ee{decompositionit}
with deterministic $\nu_t\in \bR^K$ and $R_s^t\in \bR^{K\times K}$, implying that $\I_t$  for $t \geq 2$
is an affine function of $\zeta_1,...,\zeta_t$ given by \eqref{decompositionit}; we denote it $\I_t(\zeta^t)$.
Although the reservoirs are not connected,
the inflows exhibit time and possibly space
(between regions) dependencies. Therefore,
the relevant history of the inflow process
in all regions needs to be stored in the state
vector implying that the problem cannot be solved
directly by Dynamic Programming. We allow decisions
$v_t, w_t, h_t$, and $r_t$
to depend on
$\mathcal{I}^t:=(\mathcal{I}_1,
\ldots,\mathcal{I}_t)$. Therefore, decision vector
$\bar x_t=[{h_t};v_t; {r_t};w_t]\in\bR^K\times\bR^K \times\bR^K \times \bR^K$ at stage $t$ is a function
$\bar x_t(\zeta^t)$ of disturbances $(\zeta_t)$ up to time $t$, so that system \rf{htsystem} of problem constraints
can be written as:
\begin{equation}\label{appli1}
\sum_{\tau=t-1}^t A^{t\tau}
\bar x_\tau(\zeta^\tau) \leq \bar b_t(\zeta^t) \in\bR^{m_t},\,t=1,...,N,
\end{equation}
where $A_{t,t}$ and $A_{t,t-1}$ are given matrices and the right-hand side
$\bar b_t(\zeta^t)$ is a linear function of $\zeta^t$.
\aic{Here $A^{t\tau}=0$ for $\tau<t-1$,
$$
A^{t t}=\left[\begin{array}{c|c|c|c}I_K&I_K&&\cr\hline
&-I_K&-I_K&-I_K\cr\hline
I_K&&&\cr\hline
-I_K&&&\cr\hline
&I_K&&\cr\hline
&-I_K&&\cr\hline
&&I_K&\cr\hline
&&-I_K&\cr\hline
&&&-I_K\cr
\end{array}\right]\mbox{ for }t\geq 1,\;\;
A^{t t-1}=\left[\begin{array}{c|c|c|c}-I_K&&&\cr\hline
&\ \ \ &\ \ \ &\ \ \ \cr\hline
&&&\cr\hline
&&&\cr\hline
&&&\cr\hline
&&&\cr\hline
&&&\cr\hline
&&&\cr\hline
&&&\cr
\end{array}\right] \mbox{ for }t\geq 2,
$$
and the right-hand side
$\bar b_t(\zeta^t)$ is a linear function of $\zeta^t$.
Note that for this problem we have
$n_t=4K$ and $m_t=9K$.}{}

Now, to apply the methodology of Section \ref{discrunc} it suffices
to reformulate the problem in terms of disturbances $\xi_t$ taking values in finite sets
of integers $D_t=\{1,2,\ldots,d_t\}$ of cardinality $d_t$ with disturbances $\zeta^t[\xi^t]=(\chi_{1 \xi_1},\chi_{2 \xi_2},\ldots,\chi_{t \xi_t})$ and controls $x_t(\xi^t):=\bar x_t(\chi_{1 \xi_1},\chi_{2 \xi_2},\ldots,\chi_{t \xi_t})$ thus
replacing constraints \eqref{appli1} with
\begin{equation}\label{appli2}
\sum_{\tau=t-1}^t A^{t\tau} {x}_\tau(\xi^\tau)\leq
{b}_t(\xi^t):= \bar b_t(\zeta^t[\xi^t]) \in\bR^{m_t},\;t=1,...,N,\;\forall \xi^t
\in D_1\times \ldots, D_t,
\end{equation}
which is clearly of form \rf{eq00}. Observe that Assumptions {\textbf{A.1.}} and {\textbf{A.2}} clearly hold for the reformulated system.
It is also easily seen that  Assumption {\textbf{A.3.}} holds true for
the right-hand side
${b}_t$ in \eqref{appli2} which is an additive with memory $\mu=1$ function
${b}_t(\xi^t )=\bar b_t(\zeta^t[\xi^t])$ with  $\bar b_t(\cdot)$ linear in $\zeta^t$.
Finally, let us assume that Assumption {\textbf{A.4.}}
holds for the decision rules ${x}_t$ which are restricted to be additive with memory $\mu$ functions of $\xi^t$, and let us denote
\[
u^N=\left\{u^t_{\tau\xi}\in\bR^{n_t}:1\leq\tau\leq t\leq N,\xi\in D_{\tau-\mu+1:\tau}\right\}
\]
parameter collections in the representation
\[
x_t(\xi^t)=\sum_{\tau=1}^t u_{\tau\xi_{\tau-\mu+1:\tau}}^t
\]
of candidate decision rules. Note that because the problem objective is linear in $x^N$  Assumption {\textbf{A.5.}}
obviously holds. Moreover, when the discrete distribution $P_t$ of $\zeta_t$ (and thus distribution of $\xi_t$) is known
 the objective $f(u^N)$ to be minimized in order to compute optimal constant depth decision rules is known in closed-form. Specifically, denoting
$f_t^T {x}_t(\xi^t)$ the cost per stage $t$ we have
\[
f(u^N) =\sum_{t=1}^N f_t^T \sum_{s=1}^t
 \sum_{\xi_{s-\mu+1:s} \in D_{s-\mu+1:s} } \left( \prod_{r=s-\mu+1}^s
 {P_r}(\zeta_r = \chi_{r \xi_r})\right) u_{s  \xi_{s-\mu+1:s}}^t.
\]

\subsection{Numerical experiments: comparing CDDRs and SDDP}\label{sec:num}

\par Numerical simulations described in this section utilize a MATLAB library for computing optimal Constant Depth Decision Rules for Multistage Stochastic Linear Programs with constraints of the form \eqref{eq00} satisfying the conditions in Section \ref{discrunc} with known distribution of perturbations $\xi^N$.\footnote{The functions in the library allow to load the linear program whose solutions are optimal CDDRs and solves it using Mosek \cite{mosek} solver. The library is available
at \href{https://github.com/vguigues/Constant_Depth_Decision_Rules_Library}{\tt https://github.com/vguigues/Constant\_Depth\_Decision\_Rules\_Library}. A function to run the simulations of this section %on an instance of the hydro-thermal production planning problem
is also provided, together with the implementation of SDDP for the
considered hydro-thermal application.}
\par
We compare CDDR and Stochastic Dual Dynamic Programming (SDDP, see, e.g., \cite{pereira,shapsddp,guiguescoap2013} and references therein) solutions on an instance of the hydro-thermal production planning problem with interstage dependent inflows described in the previous section.
%\footnote{SDDP solutions are implemented using the MATLAB toolbox available at {\tt\url{https://github.com/vguigues/Primal_SDDP_Library_Matlab}}.}
Parameters of the problem are initialized to mimic the mid-term Brazilian hydro-thermal problem with $K=4$ equivalent
subsystems for the country considered for instance in \cite{pythonlibrarysddp} and \cite{guiguescoap2013}.
For the sake of simplicity we considered only one equivalent thermal plant per subsystem (as described in the previous section). The state vector stores the reservoir
levels at the end of the stage along with the relevant history of inflows; the inflow model is
calibrated using available historical data with (time-dependent) model lags $(\ell_t)_{t \geq 1}$ varying between
4 and 10 months (see, e.g., \cite{shapsddp} for statistical analysis
relating the original problem and the corresponding SAA).

\subsection{Problem  without Relatively Complete Recourse}
We consider an instance of the hydro-thermal
production planning problem with $N=12$ stages (each stage representing a month)  and $d_t=10$ realizations
per stage. In each region $i$, the initial
reservoir level $v_0(i)$ is set
to half the sum of the demands for that
region over the optimization period (recall
that demand is deterministic) and we
require the final reservoir level $v_N$
to be at least $0.8v_0$ at the end of
the optimization period. Furthermore, the inflow levels are chosen in such a way that the Relatively Complete Recourse (RCR) assumption does not hold. Note that RCR is a necessary condition for application of  the standard SDDP algorithm. In order to apply SDDP, we modify the problem as follows. We  introduce slack variables
$\alpha_t$, replace constraints $v_t\geq \underline v_t$ of \rf{htsystem} with relaxed constraint
$v_t+\alpha_t\geq \underline v_t$, $\alpha_t\geq 0$, and add a penalty  $\mathrm{pen}_t^T \alpha_t$ for violation of the lower bound constraint on reservoir levels with nonnegative $\mathrm{pen}_t$. Now, we can
apply SDDP to the reformulated problem
which satisfies RCR.

We compute CDDR policies with $\mu=1,\,2$ and $3$, along with SDDP policies for four values of time-invariant penalty parameters $\mathrm{pen}_t=1,\, 10^2,\,10^3$, and $10^4$, and then simulate these
policies on 1000 scenarios of inflows. On all simulated scenarios,
CDDR policies were feasible while no SDDP policy was
feasible on all scenarios. In Table \ref{tableviolations}, for each reservoir $i$ we report the average over 1000 simulations relative violations of the lower reservoir level constraint at stage $N$ as   $\max\left(0,\frac{\underline{v}_N(i) - v_N(i)}{\underline{v}_N(i)}\right)$.
Unsurprisingly, constraint violations decrease when the penalties increase
and are very small but do not completely vanish for large values of the penalty.

For the sake of completeness we display in Table \ref{optvaltablenorcr} (without paying attention to infeasibility of the corresponding SDDP policies) simulated production costs. CPU times necessary to compute the CDDR policies and the SDDP policy corresponding to the 5\% relative suboptimality are given in Table \ref{tablecpucddrsddpnrcr}.
\begin{table}
\centering
{\small{
\begin{tabular}{|c|c|c|c|c|}
\hline
 Unit penalty $\mbox{pen}_t$   &  Reservoir 1  & Reservoir 2    &  Reservoir 3      &  Reservoir 4   \\
\hline
 1  & 0  & 0.227   &  3.8e-3     &   0   \\
\hline
 $10^2$  & 0  & 0.236   &  4e-3     &   0   \\
\hline
 $10^3$  & 0  & 2.4e-2   &  0     &   2e-4   \\
\hline
 $10^4$  & 4e-4  & 0   &  3e-4     &   0   \\
\hline
\end{tabular}
}}
\caption{Mean relative reservoir constraint violation $\max\left(0,\frac{\underline{v}_N(i) - v_N(i)}{\underline{v}_N(i)}\right)$ at  stage $N$
for reservoir $i$.}
\label{tableviolations}
\end{table}

\begin{table}
\centering
{\small{
\begin{tabular}{|c|c|c|c|c|c|c|}
\hline
 \begin{tabular}{c} CDDR\\$\mu=1$\end{tabular}    &  \begin{tabular}{c} CDDR\\$\mu=2$\end{tabular}
 & \begin{tabular}{c} CDDR\\$\mu=3$\end{tabular} &   \begin{tabular}{c}SDDP\\pen=1 \end{tabular}
 & \begin{tabular}{c}SDDP\\pen$_t$=$10^2$ \end{tabular}
 &\begin{tabular}{c}SDDP\\pen$_t$=$10^3$ \end{tabular}
 &\begin{tabular}{c}SDDP\\pen$_t$=$10^4$ \end{tabular} \\
\hline
7.81e6&6.37e6&5.87e6 &1.01e4 &1.00e4&5.2e6&5.5e6\\
\hline
\end{tabular}
}}
\caption{Simulated production costs of the CDDR and SDDP policies.}
\label{optvaltablenorcr}
\end{table}

\begin{table}
\centering
{\small{
\begin{tabular}{|c|c|c|c|c|c|c|}
\hline
\begin{tabular}{c} CDDR\\$\mu=1$\end{tabular}
&  \begin{tabular}{c} CDDR\\$\mu=2$\end{tabular}  & \begin{tabular}{c} CDDR\\$\mu=3$\end{tabular} &   \begin{tabular}{c}SDDP\\pen=1 \end{tabular}
 & \begin{tabular}{c}SDDP\\pen$_t$=$10^2$ \end{tabular}
 &\begin{tabular}{c}SDDP\\pen$_t$=$10^3$ \end{tabular}
 &\begin{tabular}{c}SDDP\\pen$_t$=$10^4$ \end{tabular} \\
\hline
4.16 & 37.8&3.1e3 &1.34e4 &4.7e4&1.03e4&4.2e3\\
\hline
\end{tabular}
}}
\caption{CPU time (in seconds) to compute CDDR and SDDP policies.}
\label{tablecpucddrsddpnrcr}
\end{table}

\subsubsection{Problem instances with Relatively Complete Recourse}
In this section we focus on the hydroelectric unit commitment problem with RCR similar to that considered in  \cite{guiguescoap2013}. When compared to the original setting of \cite{guiguescoap2013}, we consider slightly increased demands, allow less water in reservoirs at the first stage
and fix the standard deviation
of components of perturbations
$\zeta_t$ to be 0.2 (for detailed description of the set of parameters used in the simulation, see the library description).

We consider 8 instances with $N \in  \{6,12\}$ and
$d_t=d \in \{6,10,20,40\}$ (the number of realizations $d$ is the same for every stage).

For each instance, we compute optimal CDDR policies with $\mu=1,2,3$, and $4$, and run SDDP until the gap between the upper and lower bounds becomes less than 5\%.
The results of this experiment are collected in Table \ref{table1}: we report the optimal expected costs obtained using
CDDRs along with lower bound {\tt{SDDP}} LB and upper bound {\tt{SDDP}} UB at the last iteration of SDDP (results are presented for the instances in which Mosek was able to solve the corresponding LP). For all instances, we provide
in Table \ref{table2}  the CPU times along with the values of $t_{\mu=k}$---time required for SDDP to attain the optimal cost
of the CDDR policy with $\mu=k$.
%We observe that in 21 out of the 23 approximate policies
%computed (corresponding to the combinations
%of $(N,d_t,\mu)$ tested), CDDR policy is computed much
%quicker than the corresponding SDDP
%policy.
\begin{table}
\centering
{\small{
\begin{tabular}{|c|c|c|c|c|c|c|c|}
\hline
 $N$    &  $d_t$ & {\tt{SDDP}} LB   & {\tt{SDDP}} UB  &{\tt{CDDR }}$\mu=1$  & {\tt{CDDR }}$\mu=2$  & {\tt{CDDR }}$\mu=3$ &  {\tt{CDDR }}$\mu=4$     \\
\hline
6 & 6 &  2.02  &  2.10 & 2.70     & 2.30     &  2.17
&  2.07       \\
\hline
6 & 10 & 2.31  &  2.43 & 2.91     & 2.91     &  2.43
&  2.37       \\
\hline
6 & 20 & 2.04  &  2.15 & 2.71     & 2.29     &  2.17   & - \\
\hline
6 & 40 & 1.69  &  1.77 & 2.30     & 1.94   &  - &  - \\
\hline
12 & 6 & 4.92  &  5.12 & 7.01     & 5.90     &  5.41
&  -       \\
\hline
12 & 10 & 4.52  &  4.74 & 7.41     &5.90     &  5.41
&  -     \\
\hline
12& 20 & 5.22  &  5.48 & 8.17     &
6.38     & -
&  -    \\
\hline
12 & 40 &4.18  &  4.40 & 7.34    &
5.71     &  -
&  -      \\
\hline
\end{tabular}
}}
\caption{Hydro-thermal production planning example: optimal expected cost using CDDR policies
and lower and upper bound at the last iteration of SDDP.
All costs have been divided by $10^6$ to improve readability.}
\label{table1}
\end{table}

\begin{table}
\centering
{\small{
\begin{tabular}{|c|c|c|c|c|c|c|c|c|c|c|}
\hline
 $N$    &  $d_t$ & {\tt{SDDP}}  &
 \begin{tabular}{l}{\tt{CDDR }}\\$\mu=1$
\end{tabular}
 &  \begin{tabular}{l}{\tt{SDDP}}\\$t_{\mu=1}$
\end{tabular}
  &
 \begin{tabular}{l}{\tt{CDDR }}\\$\mu=2$
\end{tabular}
 &  \begin{tabular}{l}{\tt{SDDP}}\\$t_{\mu=2}$
\end{tabular}
  &
 \begin{tabular}{l}{\tt{CDDR }}\\$\mu=3$
\end{tabular}
 &  \begin{tabular}{l}{\tt{SDDP}}\\$t_{\mu=3}$
\end{tabular}
  &
 \begin{tabular}{l}{\tt{CDDR }}\\$\mu=4$
\end{tabular}
 &  \begin{tabular}{l}{\tt{SDDP}}\\$t_{\mu=4}$
\end{tabular}\\
\hline
6 & 6 & 26.1 & 0.14 & 4.4 & 1.1 & 5.5  & 4.8 &   18.4&  14.0  & 32.7 \\
\hline
6 & 10 & 37.8 & 0.6 &1.8 & 3.1 &1.8 &  22.7&  42.1& 111.2 & 156.4 \\
\hline
6 & 20 & 970  & 0.58 & 144.2 & 14.2 & 458.7 &  349.8 & 942.3 &  - & - \\
\hline
6 & 40 & 3 384 & 1.07 & 1254 & 65.5 & 2876  &  - & -  &  -  & -\\
\hline
12 & 6 & 1877 & 1.1 & 343 & 10.3 & 1115 & 144.3 & 1523  &  - & - \\
\hline
12 & 10 &1920  &  2.1 & 465  & 27.8 & 992  & 2 821.2
& 1421  &   - & -\\
\hline
12& 20 &5448  &  3.45 & 1704& 268.9 & 2358  &  - & -  & - & -  \\
\hline
12 & 40 & 15 147  & 7.5 & 2314 & 948.5 & 7836  &   - & - &  -  & - \\
\hline
\end{tabular}
}}
\caption{Hydro-thermal production planning example: CPU times (in seconds) to compute
  SDDP and CDDR policies.
Columns  {\tt{SDDP}} $t_{\mu=k}$
display the SDDP runtimes to attain the upper
bound which is equal to the optimal cost
of the CDDR policy computed with $\mu=k$.}
\label{table2}
\end{table}

\if{
\begin{table}
\centering
{\small{
\begin{tabular}{|c|c|c|c|c|}
\hline
  Depth $\mu$ &  $N$  &  $d_t$ & Nb. Variables  &Nb. constraints    \\
\hline
\hline
1 &6&6& 5 532 & 6 876  \\
\hline
2 &6&6& 22 932 & 28 476  \\
\hline
3 &6&6& 85 572 & 106 236  \\
\hline
4 &6&6& 273 492 & 339 516  \\
\hline
\hline
1 &6&10& 8 652 & 11 196  \\
\hline
2 &6&10& 58 692  &  75 996 \\
\hline
3 &6&10& 358 932 & 464 796  \\
\hline
4 &6&10& 1 860 132  & 2 408 796  \\
\hline
\hline
1 &6&20& 16 452  &  21 996  \\
\hline
2 &6&20& 220 892  &  295 596  \\
\hline
3 &6&20& 2 674 172  &  3 578 796  \\
\hline
\hline
1 &6&40& 32 052  &  43 596  \\
\hline
2 &6&40&  857 292  &  1 166 796  \\
\hline
\hline
\hline
\hline
1 &12&6& 23 592 & 119 292  \\
\hline
2 &12&6& 29 340 & 148 140  \\
\hline
3 &12&6& 589 092 & 731 340  \\
\hline
\hline
1 &12&10& 37 320 & 48 348  \\
\hline
2 &12&10& 312 540  & 404 748  \\
\hline
3 &12&10& 2 564 340 & 3 320 748  \\
\hline
\hline
1 &12&20& 71 640 & 95 868  \\
\hline
2 &12&20& 1 196 060  &   1 600 668  \\
\hline
\hline
1 &12&40& 140 280 & 190 908  \\
\hline
2 &12&40& 4 679 100  &   6 368 508  \\
\hline
\hline

\end{tabular}
}}
{\color{red}
\caption{
Number of variables and constraints of the problems solved to compute the CDDRs for the hydro-thermal application.
}
}
\label{table3}
\end{table}
}\fi

\begin{figure}
\centering
\begin{tabular}{c}
\includegraphics[scale=0.9]{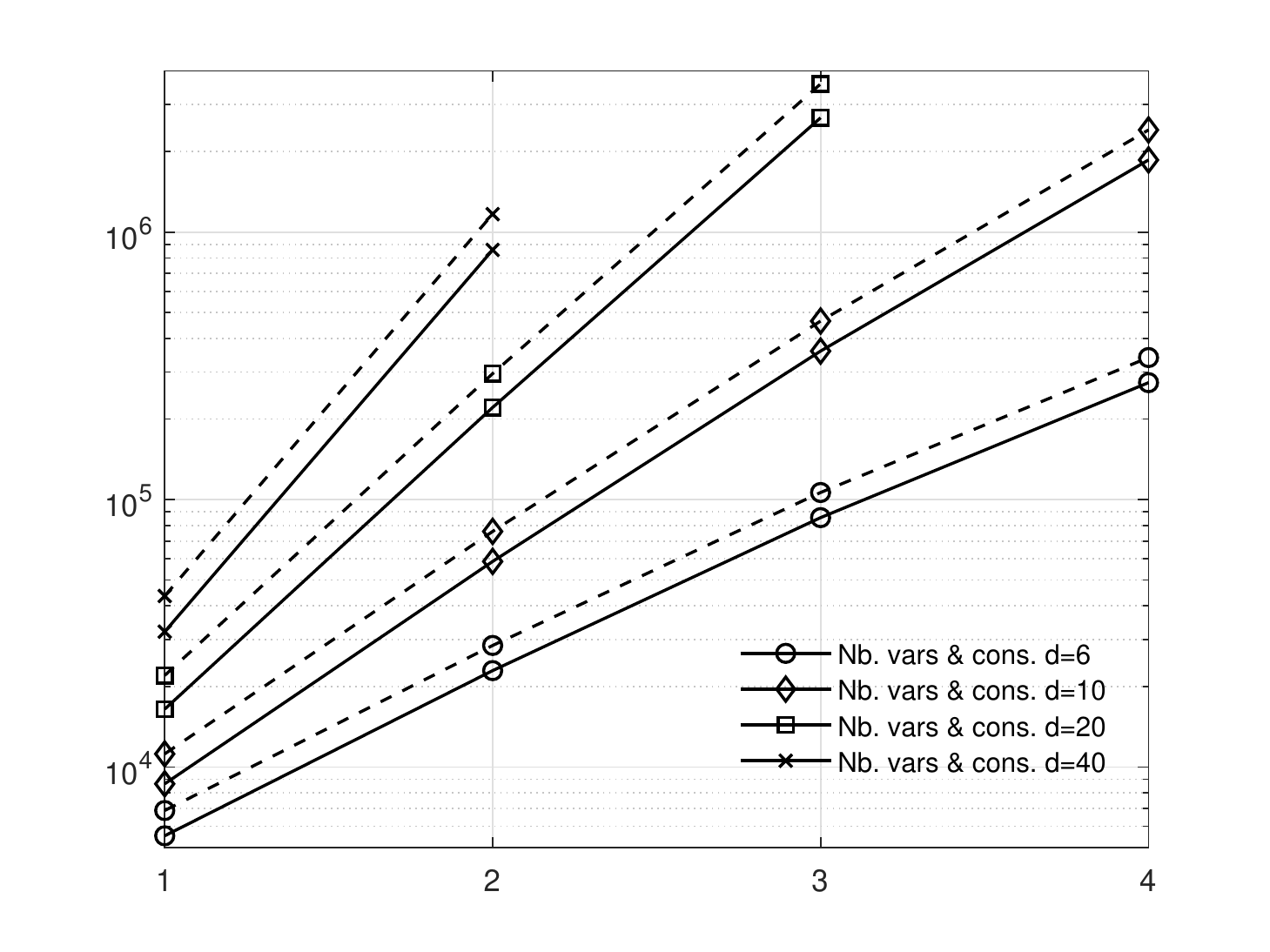}
\end{tabular}
\caption{\label{figuresize} Number of constraints (solid lines)
and variables (dashed lines) of the LPs solved to compute CDDRs for the hydro-thermal
production planning application for $d=6, 10, 20, 40$, as function of the depth $\mu$ for $N=6$.}
\end{figure}

Our experiments highlight some
strengths and drawbacks of Constant Depth Decision Rules when applied to a Multistage Stochastic Program. Compared to affine decision rules which are suboptimal in this application (their cost
cannot be less than the cost of an optimal CDDR with $\mu=1$), CDDRs are more flexible and allow
for a larger class of policies. Notice that the optimal cost of CDDRs
decreases with $\mu$ and becomes close to the lower and upper bounds computed at the last iteration of SDDP
when $\mu$ is large enough.

Stochastic Dual Dynamic Programming is now considered as tool of choice for solving Multistage Stochastic Programs. However, SDDP does not provide any guaranty  of feasibility of computed policies unless Relatively Complete Recourse condition is satisfied.
Utilizing CDDRs (as well as affine decision rules) does not require this assumption and can be efficiently used to compute feasible policies in the situation where RCR is not available. That being said, one should observe that
CDDR policies are feasible for the precise problem they are computed for, e.g., for the SAA problem on the finite
set of discrete scenarios as far as experiments in this section are concerned. Feasibility of these policies is by no means ensured
for the original hydro-thermal problem
with continuous inflows. Furthermore, the problem
to be solved to compute CDDRs can be written
for any underlying discrete process with interstage dependencies (as long as we allow memory depth
$\mu$ as large as $N-1$). CDDR also allows using continuous modeling of
disturbances with distributions supported on polytopes.
\par
As a rule of thumb, it is recommended to use CDDRs with moderate depth $\mu$ (say between $1$ and $4$) when the maximal number $d$
of possible realizations of the uncertainty for each stage is also a moderate integer and one of the above conditions applies (for instance, for problems with interstage dependent
disturbances with large lags).

The principal drawback of CDDRs is that the size of the linear optimization problem to be solved to compute the optimal rules grows exponentially with 
$\mu$ becoming prohibitive already for some toy problems considered in this section: in those examples, the values of $\mu$ for which
CDDRs could be computed for large $N$ and $d$ were too small to allow for an approximate policy of good quality.

%For instance for
%$N=12$ and $d_t=40$, we could only compute CDDRs
%for $\mu=1,2$, and the cost of the
%CDDR policy computed for $\mu=2$
%is still quite far  from SDDP cost (about 36\% %above).\\

When both SDDP and CDDR can be practiced, utilizing CDDR may be of particular interest in problems with linear constraints with random right-hand side noises satisfying autoregressive equation. When using SDDP, the size of the state vector
is proportional to the maximal lag of the autoregression in this case. On the contrary, by construction, the size of the problem to solve to compute optimal CDDRs is independent on the value of the lag. Furthermore, the numerical cost of computing SDDP policies depends on the variance of random disturbances; in general the larger the variance, the higher the cost. For CDDRs, by construction, given $d$, $\mu$,
and the number of variables and natural constraints for each
stage, the size of the optimization problem to be solved is independent of the variances
of coordinates $\zeta_t(i)$ of $\zeta_t$.

Although ``short memory'' CDDRs  (with $\mu=1$ on $\mu=2$) are suboptimal in this application, they are computed much faster than the corresponding SDDP policies, and can be used to initiate SDDP cuts.
As an illustration, we have computed, in the setting of this section, the CDDR policies with depth
$\mu=1$ for problem instances with $N=10, \,20, \,30, \,40, \,50$, and $60$ stages,
and $d_t=10$ for all stages, and ran SDDP on these instances,
stopping SDDP when the upper bound computed by SDDP attained the optimal cost by CDDR; the results are reported in Table \ref{tablecpusddpcddr}.\aic{ and show that it takes SDDP
much more time than CDDR to compute an approximate policy as good as CDDR with $\mu=1$.
This quickly computable policy (which, as we
recall, is able to compute feasible decisions
on any sample using the decision rules) could be used in combination with SDDP The corresponding CDDR policy can be used to compute
the trial points needed by SDDP for the first say hundreds forward iterations.
This CDDR policy with $\mu=1$ is at least as good as affine decision rules,
which have been used for many applications and have been shown to be optimal
for some multistage stochastic
programs, see \cite{bertsiancu2010}. We also refer
to the tests reported in Table
\ref{table2} where on most instances and for larger values
of $\mu$, SDDP needs more time
than CDDR to compute a policy
with gap not larger than the gap
of CDDR policy.}{}
\begin{table}
\centering
{\small{
\begin{tabular}{|c|c|c|c|c|c|c|}
\hline
Number of stages    &  10  & 20    & 30 & 40 & 50 & 60 \\
\hline
CDDR, $\mu=1$  & 1.6 & 7.5 & 21.4 & 60.4 & 92.1 & 157.1    \\
\hline
SDDP & 351 & 987 & 1 246 & 1 764 & 2 266 & 2 875  \\
\hline
\end{tabular}
}}
\caption{CPU time (in seconds) for CDDR with $\mu=1$ along with SDDP runtimes  needed to attain the optimal cost of the corresponding
CDDR policy.}
\label{tablecpusddpcddr}
\end{table}

As a future work, decomposition techniques could be investigated to allow using CDDRs for problems with larger values of memory depth.

\section*{Acknowledgments} Research of Vincent Guigues was supported by an
FGV grant, CNPq grants 204872/2018-9 and 401371/2014-0, FAPERJ grant E-26/201.599/2014.
Anatoli Juditsky and Arkadi Nemirovski were supported by MIAI \@ Grenoble Alpes (ANR-19-P3IA-0003), CNPq grant 401371/2014-0 and  NSF grant CCF-1523768.

\addcontentsline{toc}{section}{References}
\bibliographystyle{elsarticle-num}
\bibliography{DP_2020_0510}

\end{document}